\documentclass[10pt]{article}

\usepackage{amsmath}%

\usepackage{amsfonts}%
\usepackage{amssymb}%
\usepackage{multicol}%
\usepackage{graphicx}%
\usepackage{caption}
\usepackage{float}
\usepackage[titletoc]{appendix}

\newcommand\RR{\mathbb{R}}

\newcommand\al\alpha
\newcommand\be\beta
\newcommand\de\delta
\newcommand\ep\varepsilon
\newcommand\tha\theta
\newcommand\ka\kappa
\newcommand\la\lambda
\newcommand\om\omega

\newcommand\iy\infty
\newcommand\pa\partial

\newcommand{\hyp}[5]{\,\mbox{}_{#1}F_{#2}\!\left(\genfrac{}{}{0pt}{}{#3}{#4};#5\right)}

\numberwithin{equation}{section}

\usepackage{geometry}
\newtheorem{theorem}{Theorem}

\newtheorem{lemma}[theorem]{Lemma}

\newtheorem{Remark}[theorem]{Remark}

\begin{document}

\title{A correlation function for the classical orthogonal polynomials.}
\author{Enno Diekema \footnote{email adress: e.diekema@gmail.com}}
\maketitle

\begin{abstract}
\noindent
A correlation function of the classical orthogonal polynomials is defined and determined. The correlation function obeys a second order difference equation in two variables. The correlation function for the Gegenbauer, Chebyshev and Legendre polynomials can be written as a $_4F_3$ hypergeometric function. For the Jacobi polynomials the result is an $F_2$ Appell function. For the Generalized Laguerre polynomials the result is a confluent hypergeometric function and for the Hermite polynomials there rests only a single term.
\end{abstract}

\section{Introduction}
\setlength{\parindent}{0cm}

A well-known function in time series analysis is the correlation function of two signals. The signals can be given as two sets of data points. In that case the discrete correlation function is used. When the signals are continuous the continuous correlation function is used. In this paper the last possibility will be treated.

Given the two real functions $f(x)$ and $g(x)$ over the interval $(-\infty,\infty)$ the definition of the correlation function is
\[
R_{fg}=\int_{-\infty}^{\infty}f(x)g(x+y)dx.
\]
When $f(x)=g(x)$ this is the definition of the \text{\em autocorrelation function}. When $f(x)$ is not equal to $g(x)$ this is the definition of the \text{\em crosscorrelation function}.

In this paper the cross-correlation function is applied to the classical orthogonal  polynomials.

\

The orthogonality relation of the classical orthogonal polynomials $p_n(x)$ is given by
\[
\int_a^b p_m(x)p_n(x)w(x)dx=h_n\delta_{nm},
\]
where $\delta_{mn}$ denotes the Kronecker delta, (equal to $1$ if $m=n$ and to $0$ otherwise). An overview of the orthogonality relations is given in the following table \cite[Ch. 18]{3}.

\begin{table}[H]
	\centering
	\begin{tabular}[t]{|l|l|r|r|r|r|}
		\hline
		& & & & & \\[-2ex]		
		$p_n(x)$ &Polynomial & $a$& $b$ & $w(x)$ & $h_n$ \\[0.4ex]
		\hline
		& & & & & \\[-2ex]		$P_n(x)$ &Legendre &-1 & 1 &1 &$ \dfrac{2}{2n+1}$\\[1.4ex]
		\hline
		& & & & & \\[-2ex]
		$T_m(x)$&Chebychev of first kind&-1 & 1 &$\dfrac{1}{\sqrt{1-x^2}}$ &$\dfrac{\pi}{2} \quad n \neq 0$ \\[0.4ex]
		& & & & &$ \pi \quad n=0$\\[0.4ex]
		\hline
		& & & & & \\[-2ex]
		$U_m(x)$&Chebychev of second kind&-1 & 1 &$\sqrt{1-x^2}$ &$\dfrac{\pi}{2}$ \\[1.4ex]
		\hline
		& & & & & \\[-2ex]
		$C_n^{(\alpha)}(x)$ &Gegenbauer &-1 & 1 &$(1-x^2)^{\alpha-1/2}$ & $\dfrac{\pi 2^{1-2\alpha}\Gamma(n+2\alpha)}{n!(n+\alpha)\Gamma(\alpha)^2} \quad \alpha \neq 0 $\\[2.5ex]
		& & & & &$\dfrac{2\pi}{n^2} \quad \alpha=0$ \\[1.2ex]
		\hline
		& & & & & \\[-2ex]
		$P_n^{(\alpha,\beta)}(x)$ &Jacobi &-1 & 1 &$(1-x)^{\alpha}(1+x)^{\beta}$ & $\dfrac{2^{\alpha+\beta+1}}{2n+\alpha+\beta+1}\dfrac{\Gamma(n+\alpha+1)\Gamma(n+\beta+1)}{n!\Gamma(n+\alpha+\beta+1)} $ \\[1.8ex]
		\hline
		& & & & & \\[-2ex]
		$L_n^{(\alpha)}(x)$ &Generalized Laguerre &0 &$\infty$ &$e^{-x}x^{\alpha}$ &$\dfrac{\Gamma(\alpha+n+1)}{n!}$ \\[1.4ex]
		\hline
		& & & & & \\[-2ex]
		$H_n(x)$ &Hermite &$-\infty$ &$\infty$ &$e^{-x^2}$ & $\sqrt{\pi}\, 2^n\, n!$\\[0.4ex]
		\hline
	\end{tabular}
\end{table}

\

The orthogonal polynomials obey a three term recurrence equation. In this paper we use the following equation
\begin{equation}
p_{n+1}(x)=(B_n+A_n\, x)p_n(x)-C_n\,p_{n-1}(x),
\label{1.1}
\end{equation}
with given $p_0(x)$ and $P_1(x)$. An overview of the parameters of the recurrence equation is given in the following table

\begin{table}[H]
	\centering
	\begin{tabular}[t]{|l|l|r|r|r|r|r|}
		\hline
		& & & & & &\\[-1.2ex]		
		$p_n(x)$ &Polynomial & $A_n$&$B_n$& $C_n$ & $p_0(x)$ & $p_1(x)$ \\[2ex]
		\hline
		& & & & & &\\[-1.2ex]		
		$P_n(x)$ &Legendre &$\dfrac{2n+1}{n+1}$ & 0&$\dfrac{n}{n+1}$ &1 &$ x$\\[2ex]
		\hline
		& & & & & &\\[-1.2ex]
		$T_m(x)$&Chebychev&2 & 0& 1 &1 &x \\
		&of first kind  & & & & & \\
		\hline
		& & & & & &\\[-1.2ex]
		$U_m(x)$&Chebychev&2 & 0& 1 & 1&x \\
		&of second kind  & & & & & \\
		\hline
		& & & & & &\\[-1.2ex]
		$C_n^{(\alpha)}(x)$ &Gegenbauer &$\dfrac{2(n+\alpha)}{n+1}$ &0  &$\dfrac{n+2\alpha-1}{n+1}$ &1 &$2\alpha x$\\[2ex]
		\hline
		& & & & & &\\[-1.2ex]
		$L_n^{(\alpha)}(x)$ &Generalized &$-\dfrac{1}{n+1}$ &$\dfrac{2n+\alpha+1}{n+1}$ &$\dfrac{n+\alpha}{n+1}$ &1 &$-x+\alpha+1$\\
		&Laguerre  & & & & & \\		
		\hline
		& & & & & &\\[-1.2ex]
		$H_n(x)$ &Hermite &2 &0 &$2n$ &1 &$2x$\\[2ex]
		\hline
	\end{tabular}
\end{table}

For the Jacobi polynomials $P^{(\alpha,\beta)}_{n}(x)$ we have the following parameters
\begin{align}
&A_n=\dfrac{\Gamma(2n+\alpha+\beta+3)}{2(n+1)(n+\alpha+\beta+1)\Gamma(2n+\alpha+\beta+1)}, \nonumber \\
&B_n=\dfrac{(2n+\alpha+\beta+1)(\alpha^2-\beta^2)}{2(n+1)(n+\alpha+\beta+1)(2n+\alpha+\beta)}, \nonumber \\
&C_n=\dfrac{2(n+\alpha)(n+\beta)(2n+\alpha+\beta+2)}{2(n+1)(n+\alpha+\beta+1)(2n+\alpha+\beta)},
\label{1.2}
\end{align}
with $P_0^{(\alpha,\beta)}(x)=1$ and $P_1^{(\alpha,\beta)}(x)=\dfrac{1}{2}\big(\alpha-\beta+(\alpha+\beta+2)x\big)$.

\

For the definition of the correlation function for all orthogonal functions there are two possibilities. The first possibility starts with the correlation function without the weight function. The second possibility starts including the weight function. In this paper we use the second possibility. Then the correlation function is given by
\begin{equation}
R_{m,n}(y)=\int_a^b p_n(x)p_{n+m}(x+y)w(x)dx.
\label{1.3}
\end{equation}
For the correlation function a second order difference equation can be derived. Then there are two methods to determine the correlation function.

The first method consists of a direct calculation of the integral for the consecutive values of the parameter $m$. With a bit lucky a pattern can be seen and the general correlation function can be found. This correlation function can be checked by substitution in the second order difference equation.

The second method consists of writing the orthogonal polynomial as a series of the dependent variable. Then we can calculate the integral directly. This method generally provides a different solution than the first method.

\

This paper treats the correlation function of the classical orthogonal polynomials. These can be classified according to their orthogonality interval. Now we consider three typical cases for the orthogonality interval:

\begin{itemize}
\item  Finite orthogonality interval $[-1,1] $.

These include the Jacobi polynomials. The Gegenbauer, Chebyshev of the first and second kind and the Legendre polynomials are special cases of the Jacobi polynomials.

\item Infinite orthogonality interval $[0,\infty )$.

These include the Generalized Laguerre polynomials.

\item Infinite orthogonality interval $(-\infty ,\infty ) $.

These include the Hermite polynomials.
\end{itemize}

\

We organize this paper as follows. In Section 2 we introduce the second order recurrence equation for the correlation function. In section 3 two useful lemmas about the hypergeometric function are proved. The first is a transformation of a bounded hypergeometric function. The second lemma treats a hypergeometric function where one of the lower parameters is a negative integer. In the following sections the correlation function of  the classical orthogonal polynomials are derived.

\

\section{The second order recurrence equation of the correlation function}
To determine the recurrence equation for the correlation function, we rewrite \eqref{1.1} as
\begin{equation}
x\,p_n(x)=\dfrac{1}{A_n}p_{n+1}(x)-\dfrac{B_n}{A_n}p_n(x)+\dfrac{C_n}{A_n}p_{n-1}(x). 
\label{2.1}
\end{equation}
Multiplication both sides of this equation with the factor $x$ gives
\begin{equation}
x^2\,p_n(x)=\dfrac{1}{A_n}x\, p_{n+1}(x)-\dfrac{B_n}{A_n}x\, p_n(x)+\dfrac{C_n}{A_n}x\, p_{n-1}(x) .
\label{2.2}
\end{equation}
Application of \eqref{2.1} gives
\begin{multline}
x^2\,p_n(x)=\dfrac{1}{A_n\,A_{n+1}}p_{n+2}(x)-
\dfrac{1}{A_n}\left(\dfrac{B_{n+1}}{A_{n+1}}+\dfrac{B_n}{A_n}\right)p_{n+1}(x)+
\dfrac{1}{A_n}\left(\dfrac{C_{n+1}}{A_{n+1}}+\dfrac{B_n^2}{A_n}+\dfrac{C_n}{A_{n-1}}\right)p_n(x) \\
-\dfrac{1}{A_n}\left(\dfrac{B_n\,C_n}{A_n}+\dfrac{B_{n-1}\, C_n}{A_{n-1}}\right)p_{n-1}(x)+\dfrac{C_n\, C_{n-1}}{A_n\, A_{n-1}}p_{n-2}(x).
\label{2.3}
\end{multline}
Using \eqref{1.1} for $p_{n+m+2}(x+y)$ gives
\begin{equation}
p_{n+m+2}(x+y)=\big(B_{n+m+1}+(x+y)A_{n+m+1}\big)p_{n+m+1}(x+y)-C_{n+m+1}\,p_{n+m}(x+y).
\label{2.4}
\end{equation}
Using \eqref{1.3} for the correlation function results in
\[
R_{m+1,n+1}(y)=\int_a^b p_{n+1}(x) p_{n+m+2}(x+y)w(x)dx. 
\]
Using \eqref{1.1} for $p_{n+1}(x)$, \eqref{2.4} for $p_{n+m+2}(x)$ and \eqref{2.1} and \eqref{2.3} gives after some simplification
\begin{align}
R_{m+1,n+1}(y)&=\dfrac{A_{n+m+1}}{A_{n+1}}R_{m-1,n+2}(y)+\dfrac{A_{n+m+1}C_{n+1}}{A_{n+1}}R_{m+1,n}(y)+ \nonumber \\
&+\left(B_{n+m+1}+y\,A_{n+m+1}-\dfrac{A_{n+m+1}}{A_{n+1}}B_{n+1}\right)\,R_{m,n+1}(y)-C_{n+m+1}R_{m-1,n+1}(y).
\label{2.5}
\end{align}
This is a second order difference equation in $m$ as well as in $n$. This equation can be used to check the different solutions of the correlation function.

\

\section{Two useful lemmas}
Two lemmas can be used for the proofs of all the theorems in this paper.
\begin{lemma}
For the hypergeometric function with $n$ a non-negative integer there is the transformation
\[
_{p+1}F_{p}\left(\begin{array}{l}
-n,a_1,\dots,a_p \\
\quad \ \ \,b_1,\dots,b_p
\end{array};x\right)=
\dfrac{(a_1)_n\dots(a_p)_n}{(b_1)_n,\dots(b_p)_n}(-x)^n \ _{p+1}F_{p}\left(\begin{array}{l}
-n,1-b_1-n,\dots,1-b_p-n \\
\quad\ \ \, 1-a_1-n,\dots,1-a_p-n
\end{array};\dfrac{1}{x}\right).
\]	
\label{L3}
\end{lemma}
\underline{Proof:}\ \ The proof is straightforward by reversing the order of summation in the hypergeometric function. Starting with
\[
_{p+1}F_{p}\left(\begin{array}{l}
-n,a_1,\dots,a_p \\
\quad \ \ \,b_1,\dots,b_p
\end{array};x\right)=\sum_{i=0}^n \dfrac{(-n)_i(a_1)_i\dots(a_p)_i}{(b_1)_i\dots(b_p)i}\dfrac{1}{i!}x^i
\]
and setting $i=n-j$ results in
\[
_{p+1}F_{p}\left(\begin{array}{l}
-n,a_1,\dots,a_p \\
\quad \ \ \,b_1,\dots,b_p
\end{array};x\right)=x^n \sum_{j=0}^n \dfrac{(-n)_{n-j}(a_1)_{n-j}\dots(a_p)_{n-j}}
{(b_1)_{n-j}\dots(b_p)_{n-j}}\dfrac{(1)_j}{\Gamma(m+1-j)}\dfrac{1}{j!}\left(\dfrac{1}{x}\right)^j.
\]
Making use of 
\[
\Gamma(a-j)=(-1)^j\dfrac{\Gamma(a)}{(1-a)_j} \qquad\qquad 
(a)_{n-j}=(-1)^j\dfrac{(a)_n}{(1-a-n)_j},
\]
proofs the Lemma.
$\square $

For $p=1$ this is a well-known result.

\

\begin{lemma}
$M$ a non-negative integer. There is the transformation
\begin{multline*}
\dfrac{1}{\Gamma(-M)} \ _{p+1}F_p
\left(\begin{array}{l}
a_0,\dots,a_p \\
-M,b_2,\dots,b_p
\end{array};x\right)= \\
=\dfrac{x^{M+1}(a_0)_{M+1}\dots(a_p)_{M+1}}{\Gamma(M+2)(b_2)_{M+1}\dots(b_p)_{M+1}}
\ _{p+1}F_p
\left(\begin{array}{l}
a_0+M+1,\dots,a_p+M+1 \\
M+2,b_2+M+1,\dots,b_p+M+1
\end{array};x\right).
\end{multline*}
\label{L4}
\end{lemma}

\underline{Proof:}\ \ The proof is again straightforward. 
\begin{align*}
L&=\lim_{b_1 \rightarrow -M}\dfrac{1}{\Gamma(b_1)} \ _{p+1}F_p
\left(\begin{array}{l}
a_0,\dots,a_p \\
b_1,b_2,\dots,b_p
\end{array};x\right) \\
&=\lim_{b_1 \rightarrow -M}\dfrac{1}{\Gamma(b_1)}
\sum_{k=0}^\infty \dfrac{(a_0)_k\dots(a_p)_k}{(b_1)_k(b_2)_k\dots(b_p)_k}\dfrac{x^k}{k!} \\
&=\lim_{b_1 \rightarrow -M}\dfrac{1}{\Gamma(b_1)}
\sum_{k=0}^\infty \dfrac{\Gamma(b_1)(a_0)_k\dots(a_p)_k}{\Gamma(k+b_1)(b_2)_k\dots(b_p)_k}\dfrac{x^k}{k!} \\ 
&=\sum_{k=0}^\infty \dfrac{(a_0)_k\dots(a_p)_k}{\Gamma(k-M)(b_2)_k\dots(b_p)_k}\dfrac{x^k}{\Gamma(k+1)} \\
&=\sum_{k=0}^\infty \dfrac{(a_0)_{k+M+1}\dots(a_p)_{k+M+1}}
{\Gamma(k+1)(b_2)_{k+M+1}\dots(b_p)_{k+M+1}}
\dfrac{x^{k+M+1}}{\Gamma(k+M+2)} \\
&=\dfrac{x^{M+1}}{\Gamma(M+2)}\sum_{k=0}^\infty \dfrac{(a_0)_{k+M+1}\dots(a_p)_{k+M+1}}
{(M+2)_k(b_2)_{k+M+1}\dots(b_p)_{k+M+1}}
\dfrac{x^k}{k!}.
\end{align*}
Making use of
\[
(a)_{k+p}=(a)_p(a+p)_k
\]
gives
\[
L=\dfrac{x^{M+1}(a_0)_{M+1}\dots(a_p)_{M+1}}{\Gamma(M+2)(b_2)_{M+1}\dots(b_p)_{M+1}}
\sum_{k=0}^\infty \dfrac{(a_0+M+1)_k\dots(a_p+M+1)_k}{(M+2)_k(b_2+M+1)_k\dots(b_p+M+1)_k}
\dfrac{x^k}{k!}. 
\]
The summation can be written as a hypergeometric function. This proofs the Lemma. $\square $

\

\section{The correlation function of the Jacobi polynomials}
In this section we derive the correlation function of the Jacobi polynomials.
\begin{theorem}
The correlation function of the Jacobi polynomials $P_n^{(\alpha,\beta)}(x)$ defined as
\begin{equation}
S_J=R_{m,n}(y)=\int_{-1}^1(1-x)^\alpha (1+x)^\beta P_n^{(\alpha,\beta)}(x)P_{n+m}^{(\alpha,\beta)}(x+y)dx
\label{4.1}
\end{equation}
is given by
\begin{multline}
S_J= 2^{\alpha+\beta+1}\dfrac{\Gamma(\alpha+n+1)\Gamma(\beta+n+1)}{\Gamma(\alpha+\beta+2n+2)}
\dfrac{\Gamma(\alpha+\beta+2m+2n+1)}{\Gamma(\alpha+\beta+m+n+1)}\dfrac{1}{n!m!}\left(\dfrac{y}{2}\right)^m \\
F_2\left(\begin{array}{c}-m,\beta+n+1,-\beta-m-n \\ 
\alpha+\beta+2n+2,-\alpha-\beta-2m-2n\end{array};-\dfrac{2}{y},\dfrac{2}{y}\right)
\label{4.1a}
\end{multline}
with $y \in \RR$. 
\end{theorem}

The $F_2$ function is an $F_2$ Appell function and is defined as \cite[(16.3.2)]{3}
\[
F_2(\alpha;\beta,\beta';\gamma,\gamma';x,y)=
\sum_{m,n=0}^\infty\dfrac{(\alpha)_{m+n}(\beta)_m(\beta'')_n}{(\gamma)_m(\gamma')_n m!n!}x^my^n
\]
with $|x|+|y|<1$.

The special case of $F_2(\alpha;\beta,\beta';\gamma,\gamma';-x,x)$ gives the symmetric form
\begin{align}
F_2(\alpha;\beta,\beta';\gamma,\gamma';-x,x)
&=\sum_{k=0}^m\dfrac{(-m)_k(\beta')_k}{(\gamma')_k}\dfrac{1}{k!} x^k
\hyp32{-k,\beta,1-\gamma'-k}{\gamma,1-\beta'-k}{1} \\
&=\sum_{k=0}^m\dfrac{(-m)_k(\beta)_k}{(\gamma)_k}\dfrac{1}{k!} (-x)^k
\hyp32{-k,\beta',1-\gamma-k}{\gamma',1-\beta-k}{1}. \nonumber
\label{4.3}
\end{align}

\

\underline{Proof}:\ \ The proof of the theorem starts with writing the formula of the Jacobi polynomial as a hypergeometric summation
\[
P_n^{(\alpha,\beta)}(x)=\dfrac{(\alpha+1)_n}{n!}
\sum_{k=0}^n \dfrac{(-n)_k(\alpha+\beta+n+1)_k}{(\alpha+1)_k}\dfrac{1}{k!}
\left(\dfrac{1-x}{2}\right)^k.
\]
Substitution in \eqref{4.1} gives
\begin{multline*}
S_J=\dfrac{(a+1)_n}{n!}\dfrac{(a+1)_{n+m}}{(n+m)!}\int_{-1}^1(1-x)^\alpha(1+x)^\beta
\sum_{i=0}^n \dfrac{(-n)_i(\alpha+\beta+n+1)_i}{(\alpha+1)_i}\dfrac{1}{i!}
\left(\dfrac{1-x}{2}\right)^i \\
\sum_{j=0}^{n+m} \dfrac{(-n-m)_j(\alpha+\beta+n+m+1)_j}{(\alpha+1)_j}\dfrac{1}{j!}
\left(\dfrac{1-x-y}{2}\right)^j dx.
\end{multline*}
The summations and the integral are bounded so we can interchange the integral and the summations.
\begin{multline}
S_J=\dfrac{(a+1)_n}{n!}\dfrac{(a+1)_{n+m}}{(n+m)!}
\sum_{i=0}^n \dfrac{(-n)_i(\alpha+\beta+n+1)_i}{(\alpha+1)_i 2^i}\dfrac{1}{i!}
\sum_{j=0}^{n+m} \dfrac{(-n-m)_j(\alpha+\beta+n+m+1)_j}{(\alpha+1)_j 2^j}\dfrac{1}{j!} \\
\int_{-1}^1(1+x)^\beta(1-x)^{i+\alpha}(1-y-x)^j dx.
\label{4.2}
\end{multline}
Changing the integration variable $x$ into $2x-1$ gives an integration interval $(0,1)$.
\[
\int_{-1}^1(1+x)^\beta(1-x)^{i+\alpha}(1-y-x)^j dx=2^{\alpha+\beta+1+i}(2-y)^j
\int_0^1 x^\beta(1-x)^{\alpha+i}\left(1-\dfrac{2}{2-y}x\right)^j dx.
\]
This integral gives a hypergeometric function
\begin{multline*}
\int_{-1}^1(1+x)^\beta(1-x)^{i+\alpha}(1-y-x)^j dx= \\
2^{\alpha+\beta+1+i}(2-y)^j\dfrac{\Gamma(\alpha+1)\Gamma(\beta+1)}{\Gamma(\alpha+\beta+2)}\dfrac{(\alpha+1)_i}{(\alpha+\beta+2)_i}
\hyp21{-j,\beta+1}{\alpha+\beta+2+i}{\dfrac{2}{2-y}}.
\end{multline*}
Substitution in \eqref{4.2} and writing the hypergeometric function as a summation gives
\begin{multline*}
S_J=2^{\alpha+\beta+1}\dfrac{(a+1)_n}{n!}\dfrac{(a+1)_{n+m}}{(n+m)!}
\dfrac{\Gamma(\alpha+1)\Gamma(\beta+1)}{\Gamma(\alpha+\beta+2)}
\sum_{i=0}^n \dfrac{(-n)_i(\alpha+\beta+n+1)_i}{(\alpha+\beta+2)_i}\dfrac{1}{i!} \\
\sum_{j=0}^{n+m} \dfrac{(-n-m)_j(\alpha+\beta+n+m+1)_j}{(\alpha+1)_j}\dfrac{1}{j!} 
\left(\dfrac{2-y}{2}\right)^j\sum_{k=0}^j\dfrac{(-j)_k(\beta+1)_k}{(\alpha+\beta+2+i)_k}\dfrac{1}{k!}\left(\dfrac{2}{2-y}\right)^k.
\end{multline*}
With a property of the Pochhammer symbols we get
\[
(\alpha+\beta+2+i)_k=\dfrac{(\alpha+\beta+2)_k(\alpha+\beta+2+k)_i}{(\alpha+\beta+2)_i}.
\]
Substitution gives after interchanging the summations (and this is allowed because all the summations are bounded)
\begin{multline*}
S_J=2^{\alpha+\beta+1}\dfrac{(a+1)_n}{n!}\dfrac{(a+1)_{n+m}}{(n+m)!}
\dfrac{\Gamma(\alpha+1)\Gamma(\beta+1)}{\Gamma(\alpha+\beta+2)}
\sum_{j=0}^{n+m} \dfrac{(-n-m)_j(\alpha+\beta+n+m+1)_j}{(\alpha+1)_j}\dfrac{1}{j!} 
\left(\dfrac{2-y}{2}\right)^j \\
\sum_{k=0}^j\dfrac{(-j)_k(\beta+1)_k}{(\alpha+\beta+2)_k}\dfrac{1}{k!}\left(\dfrac{2}{2-y}\right)^k
\sum_{i=0}^n \dfrac{(-n)_i(\alpha+\beta+n+1)_i}{(\alpha+\beta+2+k)_i}\dfrac{1}{i!}.
\end{multline*}
The last summation is known
\[
\sum_{i=0}^n \dfrac{(-n)_i(\alpha+\beta+n+1)_i}{(\alpha+\beta+2+k)_i}\dfrac{1}{i!}=
\dfrac{(1-n+k)_n}{(\alpha+\beta+2+k)_n}=\dfrac{(1)_k}{\Gamma(1-n)(1-n)_k}
\dfrac{(\alpha+\beta+2)_k}{(\alpha+\beta+2+n)_k(\alpha+\beta+2)_n}.
\]
Substitution gives
\begin{multline*}
S_J=2^{\alpha+\beta+1}\dfrac{(a+1)_n}{n!}\dfrac{(a+1)_{n+m}}{(n+m)!}
\dfrac{\Gamma(\alpha+1)\Gamma(\beta+1)}{\Gamma(\alpha+\beta+n+2)}
\sum_{j=0}^{n+m} \dfrac{(-n-m)_j(\alpha+\beta+n+m+1)_j}{(\alpha+1)_j}\dfrac{1}{j!} 
\left(\dfrac{2-y}{2}\right)^j \\
\dfrac{1}{\Gamma(1-n)}\sum_{k=0}^j\dfrac{(-j)_k(\beta+1)_k(1)_k}{(1-n)_k(\alpha+\beta+n+2)_k}\dfrac{1}{k!}\left(\dfrac{2}{2-y}\right)^k.
\end{multline*}
Lemma \ref{L4} is used for the last summation. With $M+1=n$ we get
\begin{multline*}
S_J=2^{\alpha+\beta+1}\dfrac{(a+1)_{n+m}}{n!(n+m)!}
\dfrac{\Gamma(\alpha+1+n)\Gamma(\beta+1+n)}{\Gamma(\alpha+\beta+2n+2)}
\sum_{j=0}^{n+m} \dfrac{(-j)_n(-n-m)_j(\alpha+\beta+n+m+1)_j}{(\alpha+1)_j}\dfrac{1}{j!} 
\left(\dfrac{2-y}{2}\right)^j \\
\left(\dfrac{2}{2-y}\right)^n \hyp21{-j+n,\beta+1+n}{\alpha+\beta+2n+2}{\dfrac{2}{2-y}}.
\end{multline*}
The factor $(-j)_n$ only contributes if $j \geq n$. Adjustment of the summation gives
\begin{multline*}
S_J=2^{\alpha+\beta+1}\dfrac{(a+1)_{n+m}}{n!(n+m)!}
\dfrac{\Gamma(\alpha+1+n)\Gamma(\beta+1+n)}{\Gamma(\alpha+\beta+2n+2)}
\sum_{j=n}^{n+m}\dfrac{(-j)_n(-n-m)_j(\alpha+\beta+n+m+1)_j}{(\alpha+1)_j}\dfrac{1}{j!} 
\left(\dfrac{2-y}{2}\right)^j \\
\left(\dfrac{2}{2-y}\right)^n \hyp21{-j+n,\beta+1+n}{\alpha+\beta+2n+2}{\dfrac{2}{2-y}}.
\end{multline*}
Applying the new variable $k=j-n$ gives
\begin{multline*}
S_J=2^{\alpha+\beta+1}\dfrac{(a+1)_{n+m}}{n!(n+m)!}
\dfrac{\Gamma(\alpha+1+n)\Gamma(\beta+1+n)}{\Gamma(\alpha+\beta+2n+2)}
\sum_{k=0}^m\dfrac{(-n-m)_{k+n}(\alpha+\beta+n+m+1)_{k+n}}{(\alpha+1)_{k+n}}
\dfrac{(-k-n)_n}{(k+n)!} \\
\left(\dfrac{2-y}{2}\right)^k \hyp21{-k,\beta+1+n}{\alpha+\beta+2n+2}{\dfrac{2}{2-y}}.
\end{multline*}
The following properties of the Pochhammer symbols are used
\begin{align*}
&\dfrac{(-k-n)_n}{(k+n)!}=(-1)^n\dfrac{1}{k!} \\
&(-n-m)_{k+n}=(-1)^n\dfrac{(n+m)!(-m)_k}{m!}  \\
&(\alpha+\beta+n+m+1)_{k+n}=\dfrac{\Gamma(\alpha+\beta+2n+m+1)(\alpha+\beta+2n+m+1)_k}{\Gamma(\alpha+\beta+n+m+1)} \\
&(\alpha+1)_{k+n}=\dfrac{\Gamma(\alpha+1+n)(\alpha+1+n)_k}{\Gamma(\alpha+1)}.
\end{align*}
Substitution and simplifying gives
\begin{multline*}
S_J=2^{\alpha+\beta+1}\dfrac{1}{n!m!}\dfrac{\Gamma(\beta+n+1)\Gamma(\alpha+n+m+1)}{\Gamma(\alpha+\beta+2n+2)}\dfrac{\Gamma(\alpha+\beta+2n+m+1)}{\Gamma(\alpha+\beta+n+m+1)} \\
\sum_{k=0}^m\dfrac{(-m)_k(\alpha+\beta+2n+m+1)_k}{(\alpha+n+1)_k}\dfrac{1}{k!}
\left(\dfrac{2-y}{2}\right)^k \hyp21{-k,\beta+n+1}{\alpha+\beta+2n+2}{\dfrac{2}{2-y}}.
\end{multline*}
The summation with the hypergeometric function can be written as an Olsson $F_P$ function \cite[(42)]{2}. The formula is
\[
\sum_{k=0}^\infty\dfrac{(a_0)_k(a_0-c_2+1)_k}{(a_0+b_2-c_2+1)_k}\dfrac{1}{k!}
\left(\dfrac{x_2-1}{x_2}\right)
\hyp21{-k,b_1}{c_1}{\dfrac{x_1}{1-x_2}}=(x_2)^{a_0}
F_P(a_0,b_1,b_2,c_1,c_2;x_1,x_2).
\]
Application gives
\begin{multline*}
S_J=2^{\alpha+\beta+1}\dfrac{\Gamma(\beta+n+1)\Gamma(\alpha+n+m+1)}{\Gamma(\alpha+\beta+2n+2)}\dfrac{\Gamma(\alpha+\beta+2n+m+1)}{\Gamma(\alpha+\beta+n+m+1)} \\
\dfrac{1}{n!m!}\left(\dfrac{y}{2}\right)^m F_P\left(-m,\beta+n+1,-\beta-m-n,\alpha+\beta+2n+2,-\alpha-\beta-2m-2n;-\dfrac{2}{y},\dfrac{2}{y}\right).
\end{multline*}
The $F_P$ function can be written as a summation of two $F_2$ Appell functions \cite[(50) ]{2}. Because the first parameter is a negative integer, the $F_P$ function is equivalent to a single $F_2$ Appell function.
\[
F_P(-m,b_1,b_2,c_1,c_2;x_1,x_2)=\dfrac{\Gamma(-m+b_2-c_2+1)\Gamma(1-c_2)}{\Gamma(-m-c_2+1)\Gamma(b_2-c_2+1)}F_2(-m,b_1,b_2,c_1,c_2;x_1,x_2).
\]
Application proofs the theorem. $\square$

\

Note that the upper parameters and the lower parameters of the $F_2$ function meet the Saalsch\"utzian condition.

Setting $A=\alpha+n+1$ and $B=\beta+n+1$ results in the much simpler form for the $F_2$ function
\[
S_J=2^{A+B-2n-1}\dfrac{\Gamma(A)\Gamma(B)}{\Gamma(A+B)}\dfrac{\Gamma(A+B+2m-1)}{\Gamma(A+B+m-n-1)}\dfrac{1}{n!m!}\left(\dfrac{y}{2}\right)^m
F_2\left( 
\begin{array}{c}
-m,B,(1-m)-B \\
A+B,2(1-m)-(A+B)
\end{array}%
;-\dfrac{2}{y},\dfrac{2}{y}\right).
\]
Application of \eqref{4.3} gives
\begin{align*}
S_J=&2^{A+B-2n-1}\dfrac{\Gamma(A)\Gamma(B)}{\Gamma(A+B)}\dfrac{\Gamma(A+B+2m-1)}{\Gamma(A+B+m-n-1)}\dfrac{1}{n!m!}\left(\dfrac{y}{2}\right)^m \\
&\sum_{k=0}^{m-1}\dfrac{(1-m)_k(m)_k(B)_k}{(A+B)_k(2-2m-A-B)_k}\dfrac{1}{k!}\left(\dfrac{2}{y}\right)^k \hyp32{-k,m-k,B}{1-m,1-A-k}{1}
\end{align*}

\

\section{The correlation function of the Gegenbauer polynomials}
In this section we derive the correlation function of the Gegenbauer polynomials.
\begin{theorem}
The correlation function of the Gegenbauer polynomials $C^{(\alpha)}_n(x)$ defined as
\[
S_G=R_{m,n}(y)=\int_{-1}^1C^{(\alpha)}_n(x)C^{(\alpha)}_{n+m}(x+y)\left(1-x^2\right)^{\alpha-1/2}dx
\]
is given by
\begin{equation}
S_G=\dfrac{\pi\ 2^{1-2\alpha+m}\Gamma(2\alpha+n)\Gamma(\alpha+m+n)}{\Gamma(\alpha)^2\Gamma(\alpha+n+1)\Gamma(n+1)\Gamma(m+1)}\ y^m
\hyp43{-\dfrac{m}{2},\dfrac{1-m}{2},\dfrac{1-m}{2},\dfrac{2-m}{2}}{1-m,1-m-n-\alpha,n+\alpha+1}{\dfrac{4}{y^2}}
\label{5.1}
\end{equation}
with $y \in \RR$ and $m>0$.
\end{theorem}

As a first approach the definition is used. The integral is calculated for different values of $m$ and $n$.

For $m=5$ and $n=4$ Mathematica gives
\begin{align*}
S_G=R_{5,4}(y)&=\sqrt{\pi}\dfrac{\Gamma\left(\dfrac{5}{2}+\alpha\right)}{3\Gamma(5+\alpha)}\ 4\alpha^2(1+\alpha)^2(2+\alpha)(3+\alpha)(4+\alpha)y+ \\
&+\sqrt{\pi}\dfrac{\Gamma\left(\dfrac{5}{2}+\alpha\right)}{3\Gamma(5+\alpha)}\ 8\alpha^2(1+\alpha)^2(2+\alpha)(3+\alpha)(4+\alpha)(6+\alpha)(7+\alpha)y^3+ \\
&+\sqrt{\pi}\dfrac{\Gamma\left(\dfrac{5}{2}+\alpha\right)}{45\Gamma(5+\alpha)}\ 8\alpha^2(1+\alpha)^2(2+\alpha)(3+\alpha)(4+\alpha)(5+\alpha)(6+\alpha)(7+\alpha)(8+\alpha)y^5.
 \end{align*}
For $m=6$ and $n=4$ Mathematica gives
\begin{align*}
S_G=R_{6,4}(y)&=\sqrt{\pi}\dfrac{\Gamma\left(\dfrac{5}{2}+\alpha\right)}{\Gamma(5+\alpha)}\ 4\alpha^2(1+\alpha)^2(2+\alpha)(3+\alpha)(4+\alpha)(7+\alpha)y^2+ \\
&+\sqrt{\pi}\dfrac{\Gamma\left(\dfrac{5}{2}+\alpha\right)}{9\Gamma(5+\alpha)}\ 16\alpha^2(1+\alpha)^2(2+\alpha)(3+\alpha)(4+\alpha)(6+\alpha)(7+\alpha)(8+\alpha) y^4+ \\
&+\sqrt{\pi}\dfrac{\Gamma\left(\dfrac{5}{2}+\alpha\right)}{135\Gamma(5+\alpha)}\ 8\alpha^2(1+\alpha)^2(2+\alpha)(3+\alpha)(4+\alpha)(5+\alpha)(6+\alpha)(7+\alpha)(8+\alpha)(9+\alpha)y^6.
\end{align*}
Manipulation with the Gamma functions suggests
\begin{equation}
S_G=\dfrac{\pi\, 2^{1-2\alpha}\Gamma(n+2\alpha)}{n!\Gamma(\alpha)^2} \sum_{k=0}^{m-1}\dfrac{\Gamma(m+n+\alpha-k)\Gamma(m-k)}{\Gamma(n+\alpha+1+k)\Gamma(m+1-2k)\Gamma(m-2k)}\dfrac{1}{k!}(2y)^{m-2k}.
\label{5.2}
\end{equation}

\

\underline{Proof}:\ \ For the proof we use the well-known formula which indicates the relationship between the Jacobi polynomials and the Gegenbauer polynomials
\[
C_n^{(\alpha)}(x)=\dfrac{\Gamma(2\alpha+n)\Gamma\left(\alpha+\dfrac{1}{2}\right)}
{\Gamma(2\alpha)\Gamma\left(\alpha+n+\dfrac{1}{2}\right)}P_n^{(\alpha-1/2,\alpha-1/2)}(x).
\]
With this formula we get
\[
S_G=\dfrac{\Gamma(2\alpha+n)\Gamma\left(\alpha+\dfrac{1}{2}\right)}
{\Gamma(2\alpha)\Gamma\left(\alpha+n+\dfrac{1}{2}\right)}
\dfrac{\Gamma(2\alpha+n+m)\Gamma\left(\alpha+\dfrac{1}{2}\right)}
{\Gamma(2\alpha)\Gamma\left(\alpha+n+m+\dfrac{1}{2}\right)}
S_J[\alpha \rightarrow\alpha-1/2,\beta \rightarrow\alpha-1/2].
\]
Application to \eqref{4.1a} gives
\begin{equation}
S_G=\dfrac{\pi\ 2^{1-2\alpha+m}\Gamma(2\alpha+n)\Gamma(\alpha+m+n)}{\Gamma(\alpha)^2\Gamma(\alpha+n+1)\Gamma(n+1)\Gamma(m+1)}\ y^m
F_2\left(\begin{array}{c}-m,\alpha+n+\dfrac{1}{2},\dfrac{1}{2}-\alpha-m-n \\ 
2\alpha+2n+1,1-2\alpha-2m-2n\end{array};-\dfrac{2}{y},\dfrac{2}{y}\right).
\label{5.3}
\end{equation}
The $F_2$ function has a special form and for this form we use \cite[(4.4)]{4}
\begin{equation}
F_2\left(\begin{array}{c}a,b_1,b_2 \\ 
2b_1,2b_2\end{array};-x,x\right)=
\hyp43{\dfrac{a}{2},\dfrac{a+1}{2},\dfrac{b_1+b_2}{2}\dfrac{b_1+b_2+1}{2}}
{\dfrac{1+2b_1}{2},\dfrac{1+2b_2}{2},b_1+b_2}{x^2}.
\label{5.4}
\end{equation}
Application to \eqref{5.3} gives \eqref{5.1} $\square$.

The hypergeometric function is Saalsch\"utzian.

\

Equation \eqref{5.2} can be written as a hypergeometric function. We use 
\begin{equation}
\Gamma(a-k)=(-1)^k\dfrac{\Gamma(a)}{(1-a)_k}\qquad \text{and}\qquad (a)_{2k}=4^k\left(\dfrac{a}{2}\right)_k\left(\dfrac{a+1}{2}\right)_k
\label{5.5}.
\end{equation}
and the duplication formula of Legendre for the Gamma function. This results in	
\[	
S_G=\dfrac{\pi\, 2^{m+1-2\alpha}
\Gamma(n+2\alpha)\Gamma(m+n+\alpha)}{\Gamma(\alpha)^2\Gamma(n+\alpha+1)\Gamma(m+1)\Gamma(n+1)}\, y^m
\hyp43{-\dfrac{m}{2},\dfrac{1-m}{2},\dfrac{1-m}{2},\dfrac{2-m}{2}}{1-m,1-m-n-\alpha,n+\alpha+1}{\dfrac{4}{y^2}}
\]
and this is just \eqref{5.1}.

\

Another form of the correlation function is created by reversing the summation. Setting $k=m-1-r$ in \eqref{5.2} gives
\[
S_G=\dfrac{\pi\, 2^{3-m-2\alpha}\Gamma(n+2\alpha)}{n!\Gamma(\alpha)^2}y^{2-m}
\sum_{r=0}^{m-1}\dfrac{\Gamma(n+\alpha+1+r)\Gamma(r+1)\Gamma(r+1)}
{\Gamma(n+\alpha+m-r)\Gamma(3-m+2r)\Gamma(2-m+2r)\Gamma(m-r)}\dfrac{1}{r!}
(4y^2)^r.
\]
Application of the properties of the Pochhammer symbols results in
\begin{multline*}
S_G=\dfrac{\pi^2 2^{3-m-2\alpha}\Gamma(n+2\alpha)\Gamma(n+\alpha+1)}{\Gamma(n+1)\Gamma(\alpha)^2\Gamma(n+\alpha+m)\Gamma(m)\Gamma(2-m)\Gamma(3-m)}y^{2-m} \\
\sum_{r=0}^{m-1}\dfrac{(1-m)_r(1)_r(1)_r(1-n-m-\alpha)_r(n+\alpha+1)_r}
{\left(\dfrac{2-m}{2}\right)_r\left(\dfrac{3-m}{2}\right)_r\left(\dfrac{3-m}{2}\right)_r\left(\dfrac{4-m}{2}\right)_r}\dfrac{1}{r!}
\left(\dfrac{y^2}{4}\right)^r.
\end{multline*}
Writing the summation as a hypergeometric function gives
\begin{multline*}
S_G=\dfrac{\pi\, 2^{3-m-2\alpha}\Gamma(n+2\alpha)\Gamma(n+\alpha+1)}{\Gamma(n+1)\Gamma(\alpha)^2\Gamma(n+\alpha+m)\Gamma(m)\Gamma(2-m)\Gamma(3-m)}y^{2-m} \\
\hyp54{1-m,1,1,1-\alpha-m-n,n+\alpha+1}{\dfrac{2-m}{2},\dfrac{3-m}{2}\dfrac{3-m}{2},\dfrac{4-m}{2}}{\dfrac{y^2}{4}}.
\end{multline*}
In the denominator of the leading factor there are factors $\Gamma(2-m)$ and $\Gamma(3-m)$. These should be rewritten using the Legendre duplication formula. The result is
\begin{multline*}
S_G=\dfrac{\pi^2 2^{m-2\alpha}\Gamma(n+2\alpha)\Gamma(n+\alpha+1)}
{\Gamma(n+1)\Gamma(\alpha)^2\Gamma(n+\alpha+m)\Gamma(m)\Gamma\left(\dfrac{2-m}{2}\right)\Gamma\left(\dfrac{3-m}{2}\right)\Gamma\left(\dfrac{3-m}{2}\right)\Gamma\left(\dfrac{4-m}{2}\right) }y^{2-m} \\
\hyp54{1-m,1,1,1-\alpha-m-n,n+\alpha+1}{\dfrac{2-m}{2},\dfrac{3-m}{2}\dfrac{3-m}{2},\dfrac{4-m}{2}}{\dfrac{y^2}{4}}.
\end{multline*}
This form is given by Mathematica summing \eqref{5.2}.

\

Using Lemma \ref{L4} with $m$=even and $M=\dfrac{m-2}{2}$ gives
\[
S_G=\dfrac{\pi\, 2^{-2\alpha}
\Gamma(n+2\alpha)}{\Gamma(\alpha)^2\Gamma(n+1)}m(m+2n+2\alpha)y^2
\hyp43{\dfrac{2-m}{2},\dfrac{2+m}{2},\dfrac{2-m}{2}-n-\alpha,\dfrac{2+m}{2}+n+\alpha}{\dfrac{3}{2},\dfrac{3}{2},2}{\dfrac{y^2}{4}}.
\]

Using Lemma \ref{L4} with $m$=odd and $M=\dfrac{m-3}{2}$ gives
\[
S_G=\dfrac{\pi\, 2^{2-2\alpha}\Gamma(n+2\alpha)}{\Gamma(\alpha)^2\Gamma(n+1)}y
\hyp43{\dfrac{1-m}{2},\dfrac{1+m}{2},\dfrac{1-m}{2}-n-\alpha,\dfrac{1+m}{2}+n+\alpha}{\dfrac{1}{2},1,\dfrac{3}{2}}{\dfrac{y^2}{4}}.
\]

\

\section{The correlation function of the Chebyshev polynomials of the first kind}

In this section we derive the correlation function of the Chebyshev polynomials of the first kind
.
\begin{theorem}
The correlation function of the Chebyshev polynomials of the first kind $T_n(x)$ defined as
\[
S_T=R_{m,n}(y)=\int_{-1}^1T_n(x)T_{n+m}(x+y)\dfrac{1}{\sqrt{1-x^2}}dx
\]
is given by
\begin{equation}
S_T=\dfrac{\pi\Gamma(m+n+1)}{\Gamma(m+1)\Gamma(n+1)}2^{m-1}y^m
\hyp43{-\dfrac{m}{2},\dfrac{1-m}{2},\dfrac{1-m}{2},\dfrac{2-m}{2}}{1-m,1-m-n,n+1}{\dfrac{4}{y^2}}.
\label{6.1}
\end{equation}
with $y \in \RR$.	
\end{theorem}

\

As a first approach the definition is used. The integral is calculated for different values of $m$ and $n$.

For $m=8$ and $n=4$ Mathematica gives
\[
S_T=R_{5,4}(y)=\pi\big(384 y^2+20160 y^4+96768 y^6+63360 y^8\big).
\]
For $m=9$ and $n=4$ Mathematica gives
\[
S_T=R_{6,4}(y)=\pi\big(13 y+6240 y^3+131040 y^5+384384 y^7+183040 y^9\big).
\]
Manipulation with the Gamma functions suggests
\begin{equation}
S_T= \pi \left(\dfrac{n+m}{2}\right)\sum_{k=0}^{m-1}\dfrac{\Gamma(m+n-k)\Gamma(m-k)}
{\Gamma(n+1+k)\Gamma(m+1-2k)\Gamma(m-2k)}\dfrac{1}{k!}(2y)^{m-2k}.
\label{6.2}
\end{equation}

\

\underline{Proof}:\ \ For the proof we use the well-known formula which indicates the relationship between the Gegenbauer polynomials and the Chebyshev polynomials of the first kind.
\[
T_n(x)=\dfrac{P_n^{(-1/2,-1/2)}(x)}{P_n^{(-1/2,-1/2)}(1)}.
\]
For the denominator there is
\[
P_n^{(-1/2,-1/2)}(1)=\dfrac{\Gamma(2n+1)}{2^{2n}\Gamma(n+1)^2}
\]
so we get
\[
T_n(x)=\dfrac{2^{2n}\Gamma(n+1)^2}{\Gamma(2n+1)}P_n^{(-1/2,-1/2)}(x).
\]
With this formula we get
\begin{equation}
S_T=\dfrac{2^{4n+2m}\Gamma(n+1)^2\Gamma(n+m+1)^2}{\Gamma(2n+1)\Gamma(2n+2m+1)}
S_J[\alpha \rightarrow -1/2,\ \beta \rightarrow -1/2].
\label{6.3}
\end{equation}
Application of \eqref{4.1a} gives
\begin{align*}
S_J[\alpha \rightarrow -1/2,\ \beta \rightarrow -1/2]
&=\dfrac{\Gamma\left(n+\dfrac{1}{2}\right)^2}{\Gamma(2n+1)}
\dfrac{\Gamma(2m+2n)}{\Gamma(m+n)}\dfrac{1}{n!m!}\left(\dfrac{y}{2}\right)^m 
F_2\left(\begin{array}{c}-m,n+\dfrac{1}{2},\dfrac{1}{2}-m-n \\ 
2n+1,1-2m-2n\end{array};-\dfrac{2}{y},\dfrac{2}{y}\right) \\
&=\dfrac{\Gamma\left(n+\dfrac{1}{2}\right)\Gamma\left(m+n+\dfrac{1}{2}\right)2^{2m-1}}{\Gamma(n+1)^2\Gamma(m+1)}
\left(\dfrac{y}{2}\right)^m 
F_2\left(\begin{array}{c}-m,n+\dfrac{1}{2},\dfrac{1}{2}-m-n \\ 
2n+1,1-2m-2n\end{array};-\dfrac{2}{y},\dfrac{2}{y}\right) \\
&=\dfrac{\pi\Gamma(2n+1)\Gamma(2m+2n+1)}{\Gamma(n+1)^3\Gamma(m+1)\Gamma(m+n+1)
2^{4n+1}}
\left(\dfrac{y}{2}\right)^m 
F_2\left(\begin{array}{c}-m,n+\dfrac{1}{2},\dfrac{1}{2}-m-n \\ 
	2n+1,1-2m-2n\end{array};-\dfrac{2}{y},\dfrac{2}{y}\right).
\end{align*}
Substitution in \eqref{6.3} results in
\[
S_T=\dfrac{\pi\Gamma(m+n+1)}{\Gamma(n+1)\Gamma(m+1)}\ 2^{m-1}y^m
F_2\left(\begin{array}{c}-m,n+\dfrac{1}{2},\dfrac{1}{2}-m-n \\ 
2n+1,1-2m-2n\end{array};-\dfrac{2}{y},\dfrac{2}{y}\right).
\]
Application of \ref{5.4} completes the proof $\square$.

The hypergeometric function is Saalsch\"utzian.

\

Equation \eqref{6.2} can be written as a hypergeometric function. We use \eqref{5.5}
and the duplication formula of Legendre for the Gamma function. This results in	
\[	
S_T=\dfrac{\pi\Gamma(m+n+1)}{\Gamma(m+1)\Gamma(n+1)}2^{m-1}y^m
\hyp43{-\dfrac{m}{2},\dfrac{1-m}{2},\dfrac{1-m}{2},\dfrac{2-m}{2}}{1-m,1-m-n,n+1}{\dfrac{4}{y^2}}
\]
and this is just \eqref{6.1}.

\

Reversing the order of summation in \eqref{6.2} results in 
\[
S_T=\pi\ m(m+n)(m+2n)\dfrac{y^2}{4} \hyp43{\dfrac{2-m}{2},\dfrac{2+m}{2},\dfrac{2-m}{2}-n,\dfrac{2+m}{2}+n}{\dfrac{3}{2},\dfrac{3}{2},2}{\dfrac{y^2}{4}}
\]
for $m$ even and for $m$ odd
\[
S_T=\pi\ (m+n)y
\hyp43{\dfrac{1-m}{2},\dfrac{1+m}{2},\dfrac{1-m}{2}-n,\dfrac{1+m}{2}+n}{\dfrac{1}{2},1,\dfrac{3}{2}}{\dfrac{y^2}{4}}
\]
with $y \in \RR$. Both hypergeometric functions are Saalsch\"utzian.

\

Mathematica gives the following result for summing \eqref{6.2}
\[
S_T=\dfrac{\pi\, 2^{1-m}y^{2-m}(m+n)\Gamma(n+1)}{\Gamma(2-m)\Gamma(3-m)\Gamma(m)\Gamma(m+n)}
\hyp54{1-m,1,1,1-m-n,1+n}{\dfrac{2-m}{2},\dfrac{3-m}{2},\dfrac{3-m}{2},\dfrac{4-m}{2}}{\dfrac{y^2}{4}}.
\]
This can be proven in the same way as in the previous section.

\

\section{The correlation function of the Chebyshev polynomials of the se\-cond kind}

In this section we derive the correlation function of the Chebyshev polynomials of the second kind
\begin{theorem}
The correlation function of the Chebyshev polynomials of the second kind $U_n(x)$ defined as
\[
S_U=R_{m,n}(y)=\int_{-1}^1U_n(x)U_{n+m}(x+y)\sqrt{1-x^2}dx
\]
is given by
\begin{equation}
S_U=\dfrac{\pi\Gamma(m+n+1)}{\Gamma(m+1)\Gamma(n+1)}2^{m-1}y^m
\hyp43{-\dfrac{m}{2},\dfrac{1-m}{2},\dfrac{1-m}{2},\dfrac{2-m}{2}}{1-m,1-m-n,n+1}{\dfrac{4}{y^2}}
\label{7.1}
\end{equation}
with $y \in \RR$.	
\end{theorem}

\

As a first approach the definition is used. The integral is calculated for different values of $m$ and $n$.

For $m=8$ and $n=4$ Mathematica gives
\[
S_U=R_{8,4}(y)=\pi\big(180 y^2+12000 y^4+73920 y^6+63360 y^8\big).
\]
For $m=9$ and $n=4$ Mathematica gives
\[
S_U=R_{9,4}(y)=\pi\big(5 y+3000 y^3+79200 y^5+295680 y^7+183040 y^9\big).
\]
Manipulation with Gamma functions suggests
\begin{equation}
S_U= \pi \left(\dfrac{n+1}{2}\right)\sum_{k=0}^{m-1}\dfrac{\Gamma(m+n+1-k)\Gamma(m-k)}
{\Gamma(n+2+k)\Gamma(m+1-2k)\Gamma(m-2k)}\dfrac{1}{k!}(2y)^{m-2k}.
\label{7.2}
\end{equation}

\

\underline{Proof}:\ \ For the proof we use the well-known formula which indicates the relationship between the Gegenbauer polynomials and the Chebyshev polynomials of the second kind.
\[
U_n(x)=C_n^{(1)}(x).
\]
With this formula we get
\begin{equation}
S_U=S_C[\alpha \rightarrow 1].
\label{7.3}
\end{equation}
Application of \eqref{5.1} gives
\begin{equation}
S_U=\dfrac{\pi\ 2^{m-1}\Gamma(n+2)\Gamma(m+n+1 )}{\Gamma(\alpha)^2\Gamma(\alpha+n+1)\Gamma(n+1)\Gamma(m+1)}\ y^m
\hyp43{-\dfrac{m}{2},\dfrac{1-m}{2},\dfrac{1-m}{2},\dfrac{2-m}{2}}{1-m,1-m-n-\alpha,n+\alpha+1}{\dfrac{4}{y^2}}.
\end{equation}

This completes the proof $\square$.

The hypergeometric function is Saalsch\"utzian.

\

Mathematica gives the following result for summing \eqref{7.2}
\[
S_U=\dfrac{\pi\, 2^{1-m}y^{2-m}(n+1)\Gamma(n+2)}{\Gamma(2-m)\Gamma(3-m)\Gamma(m)\Gamma(m+n+1)}
\hyp54{1-m,1,1,-m-n,2+n}{\dfrac{2-m}{2},\dfrac{3-m}{2},\dfrac{3-m}{2},\dfrac{4-m}{2}}{\dfrac{y^2}{4}}.
\]
This can be proven in the same way as in a previous section.

\

\section{The correlation function of the Legendre polynomials}

In this section we derive the correlation function of the Legendre polynomials
\begin{theorem}
The correlation function of the Legendre polynomials $P_n(x)$ defined as
\[
S_P=R_{m,n}(y)=\int_{-1}^1P_n(x)P_{n+m}(x+y)dx
\]
is given by
\begin{equation}
S_P=\dfrac{\Gamma\left(m+n+\dfrac{1}{2}\right)}{\Gamma(m+1)\Gamma\left(n+\dfrac{3}{2}\right)}(2y)^m\hyp43{-\dfrac{m}{2},\dfrac{1-m}{2},\dfrac{1-m}{2},\dfrac{2-m}{2}}{1-m,\dfrac{1}{2}-m-n,n+\dfrac{3}{2}}{\dfrac{4}{y^2}}
\label{8.1}
\end{equation}
with $y \in \RR$.	
\end{theorem}

\

As a first approach the definition is used. The integral is calculated for different values of $m$ and $n$.

For $m=8$ and $n=4$ Mathematica gives
\[
S_P=R_{8,4}(y)=68y^2+\dfrac{8075}{2}y^4+\dfrac{88179}{4}y^6+\dfrac{1062347}{64}y^8.
\]
For $m=9$ and $n=4$ Mathematica gives
\[
S_P=R_{9,4}(y)=2y+\dfrac{3230}{3}y^3+\dfrac{101745}{4}y^5+\dfrac{676039}{8}y^7+
\dfrac{26558675}{64}y^9.
\]
Manipulation with Gamma functions suggests
\begin{equation}
S_P= \sum_{k=0}^{m-1}\dfrac{\Gamma\left(n+m+\dfrac{1}{2}-k\right)\Gamma(m-k)}
{\Gamma\left(n+\dfrac{3}{2}+k\right)\Gamma(m+1-2k)\Gamma(m-2k)}\dfrac{1}{k!}(2y)^{m-2k}.
\label{8.2}
\end{equation}

\

\underline{Proof}:\ \ For the proof we use the well-known formula which indicates the relationship between the Gegenbauer polynomials and the Legendre polynomials.
\[
P_n(x)=C_n^{(1/2)}(x)
\]
With this formula we get
\begin{equation}
S_P=S_C\left[\alpha \rightarrow \dfrac{1}{2}\right].
\label{8.3}
\end{equation}
Application of \eqref{5.1} gives
\begin{equation}
S_P=\dfrac{\Gamma\left(m+n+\dfrac{1}{2}\right)}{\Gamma(m+1)\Gamma\left(n+\dfrac{3}{2}\right)}(2y)^m\hyp43{-\dfrac{m}{2},\dfrac{1-m}{2},\dfrac{1-m}{2},\dfrac{2-m}{2}}{1-m,\dfrac{1}{2}-m-n,n+\dfrac{3}{2}}{\dfrac{4}{y^2}}.
\label{8.4}
\end{equation}

This completes the proof $\square$.

The hypergeometric function is Saalsch\"utzian.

\

There appear to be several equivalent forms of the correlation function 
\begin{align}
S_P
&=\dfrac{\Gamma\left(m+n+\dfrac{1}{2}\right)}{\Gamma(m+1)\Gamma\left(n+\dfrac{3}{2}\right)}(2y)^m\hyp43{-\dfrac{m}{2},\dfrac{1-m}{2},\dfrac{1-m}{2},\dfrac{2-m}{2}}{1-m,\dfrac{1}{2}-m-n,n+\dfrac{3}{2}}{\dfrac{4}{y^2}} \qquad n \geq 0 \nonumber \\
&=2 \dfrac{(2n+2)_{2m-1}}{(n)_m}\dfrac{1}{m!}\left(\dfrac{y}{2}\right)^m
\hyp43{-\dfrac{m}{2},\dfrac{1-m}{2},\dfrac{1-m}{2},\dfrac{2-m}{2}}{1-m,\dfrac{1}{2}-m-n,n+\dfrac{3}{2}}{\dfrac{4}{y^2}} \nonumber \\
&=\sum_{k=0}^{m-1}\binom{n+\dfrac{1}{2}+m-1+k}{n+\dfrac{1}{2}+k}\binom{m-1-k}{k}
\dfrac{1}{(m-2k)}(2y)^{m-2k}  \nonumber \\
&=\sum_{k=0}^{m-1}\dfrac{\Gamma\left(n+m+\dfrac{1}{2}-k\right)\Gamma(m-k)}
{\Gamma\left(n+\dfrac{3}{2}+k\right)\Gamma(m+1-2k)\Gamma(m-2k)}\dfrac{1}{k!}(2y)^{m-2k}.
\label{8.5}
\end{align}
If one of these forms is proven, it is easy to deduce the other forms from the proven form. 

Reversing the order of summation in \eqref{8.2} results in 
\[
S_P=\dfrac{m}{2}(m+2n+1)\, y^2\hyp43{\dfrac{2-m}{2},\dfrac{2+m}{2},\dfrac{3+m}{2}+n,\dfrac{1-m}{2}-n}{\dfrac{3}{2},\dfrac{3}{2},2}{\dfrac{y^2}{4}}
\]
for $m$ even and for $m$ odd
\[
S_P=2y\hyp43{\dfrac{1-m}{2},\dfrac{1+m}{2},\dfrac{m}{2}+n+1,-\dfrac{m}{2}-n}{1,\dfrac{1}{2},\dfrac{3}{2}}{\dfrac{y^2}{4}}
\]
with $y \in \RR$. Both hypergeometric functions are Saalsch\"utzian.

\

Mathematica gives the following result for summing \eqref{8.2}
\[
S_P=\dfrac{(2y)^{2-m}\Gamma\left(n+\dfrac{3}{2}\right)}
{\Gamma(2-m)\Gamma(3-m)\Gamma(m)\Gamma\left(m+n+\dfrac{1}{2}\right)}
\hyp54{1-m,1,1,\dfrac{1}{2}-m-n,\dfrac{3}{2}+n}{\dfrac{2-m}{2},\dfrac{3-m}{2},\dfrac{3-m}{2},\dfrac{4-m}{2}}{\dfrac{y^2}{4}}.
\]
This can be proven in the same way as in the previous sections.

\

\section{The correlation function of the Generalized Laguerre polynomials}

In this section we derive the correlation function of the Generalized Laguerre polynomials.
\begin{theorem}
The correlation function of the Generalized Laguerre polynomials $L^{(\alpha)}_n(x)$ defined as
\[
S_L=R_{m,n}(y)=\int_0^\infty L^{(\alpha)}_n(x)L^{(\alpha)}_{n+m}(x+y)e^{-x}x^\alpha dx
\]
is given by
\begin{equation}
S_L=-\dfrac{\Gamma(n+1+\alpha)}{n!}y\, \hyp11{1-m}{2}{y}
\label{9.1}
\end{equation}
with $y \in \RR$.	
\end{theorem}
Note that the hypergeometric function can be written as a Laguerre polynomial.

\

As a first approach the definition is used. The integral is calculated for different values of $m$ and $n$.

For $m=7$ and $n=4$ Mathematica gives
\[
S_L=R_{7,4}(y)=-\Gamma(5+\alpha)\left(-\dfrac{1}{24}y+\dfrac{1}{8}y^2-\dfrac{5}{48}y^3+
\dfrac{5}{144}y^4-\dfrac{1}{192}y^5+\dfrac{1}{2880}y^6-\dfrac{1}{120960}y^7 \right).
\]

Manipulation with Gamma functions suggests
\begin{equation}
S_L=\dfrac{\Gamma(n+1+\alpha)}{n!}\sum_{k=0}^m{\binom{m-1}{k-1}\dfrac{(-y)^k}{k!}}=
-\dfrac{\Gamma(n+1+\alpha)}{n!}\sum_{k=0}^{m-1}\dfrac{(1-m)_k}{(2)_k}\dfrac{y^{k+1}}{k!}.
\label{9.2}
\end{equation}

\

\underline{Proof}:\ \ For the proof of \eqref{9.2} we use the difference equation \eqref{2.5}. Substitution of the values of $A_n$, $B_n$ and $C_n$ in the difference equation gives
\begin{multline*}
(n+m+2)R_{m+1,n+1}(y)-(n+2)R_{m-1,n+2}(y)-(n+\alpha+1)R_{m+1,n}(y)- \\
-(2m-y)R_{m,n+1}(y)+(n+m+\alpha+1)R_{m-1,n+1}(y)=0.
\end{multline*}
Substitution of the right hand-side of \eqref{9.2} gives after some simplification
\[
(m+1)\sum_{k=0}^m\dfrac{(-m)_k}{(2)_k}\dfrac{y^k}{k!}-
2m\sum_{k=0}^{m-1}\dfrac{(1-m)_k}{(2)_k}\dfrac{y^k}{k!}+
\sum_{k=0}^{m-1}\dfrac{(1-m)_k}{(2)_k}\dfrac{y^{k+1}}{k!}+
(m-1)\sum_{k=0}^{m-2}\dfrac{(2-m)_k}{(2)_k}\dfrac{y^k}{k!}=0.
\]
Comparing the coefficients of the powers of $y$ gives 4 cases.
\vspace{2mm}

Case 1. Terms with $y^0$.
The resulting equation becomes:  $(m+1)-2m+(m+1)=0$.
\vspace{2mm}

Case 2. Terms with $y^r$ with $1 \leq r \leq m-2$. The resulting equations becomes
\[
(m+1)\dfrac{(-m)_r}{(2)_r}-2m\dfrac{(1-m)_r}{(2)_r}+\dfrac{(1-m)_{r-1}}{(2)_{r-1}}+
(m-1)\dfrac{(2-m)_r}{(2)_r}=0.
\]
After some manipulation of the Pochhammer symbols this gives an identity.
\vspace{2mm}

Case 3. Terms with $y^{m-1}$. The resulting equations becomes
\[
(m+1)(-m)_{m-1}-2m(1-m)_{m-1}+(m-1)(1-m)_{m-2}\dfrac{(2)_{m-1}}{(2)_{m-2}}=0.
\]
After some manipulation of the Pochhammer symbols this gives an identity.
\vspace{2mm}

Case 4.	Terms with $y^m$. The resulting equations becomes
\[
(m+1)\dfrac{(-m)_m}{(2)_m}+m\dfrac{(1-m)_{m-1}}{(2)_{m-1}}=0.
\]
After some manipulation of the Pochhammer symbols this gives an identity.

The right hand-side of equation \eqref{9.2} is another way of writing the right hand-side of equation \eqref{9.1}.

\

Because the difference equation is  linear we had to check an initial value.
$n=0,m=0$ gives
\[
R_{0,0}(y)=\int_0^\infty L^{(\alpha)}_0(x)L^{(\alpha)}_{0}(x+y)e^{-x}x^\alpha dx.
\]
Because $ L^{(\alpha)}_0(x)=1$ there follows
\[
R_{0,0}(y)=\int_0^\infty e^{-x}x^\alpha dx=\Gamma(\alpha+1).
\]
This is the same as the left hand-side of equation \eqref{9.2}. This proves the theorem. $\square$

\

For the Generalized Laguerre polynomials an addition formula is known \cite[page 85]{1}:
\[
L_n^{\alpha+\beta+1}(x+y)=\sum_{k=0}^nL_k^{(\alpha)}(x)L_{n-k}^{(\beta)}(y)
\]
This formula can also be used to prove the Theorem.

\

\section{The correlation function of the Hermite polynomials}

In this section we derive the correlation function of the Hermite polynomials.
\begin{theorem}
The correlation function of the Hermite polynomials $H^{(\alpha)}_n(x)$ defined as
\[
S_H=R_{m,n}(y)=\int_{-\infty}^\infty H_n(x)H_{n+m}(x+y)e^{-x^2}dx
\]
is given by
\begin{equation}
S_H=2^{n+m}\sqrt{\pi}(m+1)_n y^m
\label{10.1}
\end{equation}
with $y \in \RR$.	
\end{theorem}

\

The integral can be computed with Mathematica. But a check with the difference equation is needed.

\

\underline{Proof}:\ \ For the proof of \eqref{10.1} we use the difference equation \eqref{2.5}. Using the values of $A_n$, $B_n$ and $C_n$ the difference equation becomes
\[
R_{m+1,n+1}(y)-R_{m-1,n+2}(y)-2(n+1)R_{m+1,n}(y)-2yR_{m,n+1}(y)+2(n+m+1)R_{m-1,n+1}(y)=0.
\]
Substitution of the right hand-side of \eqref{10.1} gives after some simplification
\[
2\big((m+2)_{n+1}-(n+1)(m+2)_n-(m+1)_{n+1}\big)y^2+\big((n+m+1)(m)_{n+1}-(m)_{n+2}\big)=0.
\]
Comparing the coefficients of the powers of $y$ gives 2 cases.
\vspace{2mm}

Case 1. Terms with $y^2$. The resulting equation becomes  $(m+2)_{n+1}-(n+1)(m+2)_n-(m+1)_{n+1}=0$. Application of the properties of the Pochhammer symbols gives an identity.
\vspace{2mm}

Case 2. Terms with $y^0$ The resulting equations becomes
$(n+m+1)(m)_{n+1}-(m)_{n+2}=0$. Application of the properties of the Pochhammer symbols gives also an identity.
\vspace{2mm}

Because the difference equation is  linear we had to check an initial value.
$n=0,m=0$ gives
\[
R_{0,0}(y)=\int_{-\infty}^\infty H_0(x)H_{0}(x+y)e^{-x^2}dx.
\]
Because $ H_0(x)=1$ there follows
\[
R_{0,0}(y)=\int_{-\infty}^\infty e^{-x^2}dx=\sqrt{\pi}.
\]
This is the same as the right hand-side of equation \eqref{10.1} with $m=0$ and $n=0$. This proves the theorem. $\square$

\
For the Hermite polynomials an addition formula is known \cite[page 82]{1}:
\[
H_n(x+y)=\sum_{k=0}^n\binom{n}{m}H_k(x)(2y)^{n-m}.
\]
This formula can also be used to prove the Theorem.

\

\textbf{Acknowledgements}

The author is very grateful to Dr. T.H. Koornwinder and Dr. N.M. Temme for their careful and thorough reading of this paper.

\

\end{document}